\documentclass[12pt]{amsart}

\usepackage{tikz}
\usetikzlibrary{calc,through,backgrounds}
\usepackage{multicol}
\usepackage{bm,amssymb,mathtools,amsthm,amsfonts,amsmath,latexsym,cite,color, psfrag,graphicx,ifpdf,ytableau,pgfkeys,pgfopts,xcolor,extarrows}

\newtheorem{theorem}{Theorem}[section]
\newtheorem*{theorem*}{Theorem}

\newtheorem{lemma}[theorem]{Lemma}
\newtheorem{proposition}[theorem]{Proposition}
\newtheorem{corollary}[theorem]{Corollary}

\newtheorem*{conjecture*}{Conjecture}

\newtheorem{example}[theorem]{Example}
\newtheorem{remark}[theorem]{Remark}
\newtheorem{definition}[theorem]{Definition}

\allowdisplaybreaks

\makeatletter \@addtoreset{equation}{section}
\@addtoreset{theo}{}\makeatother

\usepackage{hyperref}

\renewcommand{\P}{\mathbb{P}}

\renewcommand{\hat}[1]{\widehat{#1}}

\begin{document}

\title[QCOHA of the $d$-loop quiver revisited]{Quantized cohomological Hall algebra of the $d$-loop quiver revisited}\thanks{Supported by the National Science Foundation of China (Grant No. 11821001, 11971326, 12071320)
and  the Sichuan Science and Technology Program (Grant No. 2020YJ0006).}

\author{Neil J.Y. Fan,  Changjian Fu and Liangang Peng}
\address{Department of Mathematics\\SiChuan University\\610064 Chengdu\\P.R.China}
\email{fan@scu.edu.cn(N.J.Y. Fan)}
\email{changjianfu@scu.edu.cn(C. Fu)}
\email{penglg@scu.edu.cn (L. Peng)}

\date{}

\maketitle %不能删

% 文章的正文自此开始。注意：
% 1. 这里有用 pstricks 软件包画的交换图和其他图标的范例，原指令和编译后的效果可对照参考（可用此方法画图，也可不用）；
%    如果选用其他（非LaTeX主程序内含的）软件包，在本文件的引言部分应该用\usepackage{}引用。
% 2. 这里还有命题、公式、参考文献的标签和交叉引用的例子，也供参考（论文较短时，也许就不需要这么麻烦了）。
% 3. 这里还有为编辑全书索引留下的标签，请把要做索引的关键词在正常行文之后插入\index{重复一下关键词}，便于后期处理。

\begin{abstract}Let $\Lambda$ be the set of  partitions of length $\geq 0$.
We introduce an $\mathbb{N}$-graded algebra $\mathbb{A}_q^d(\Lambda)$ associated to $\Lambda$, which can be viewed as a quantization of the algebra of partitions defined by Reineke. The multiplication of $\mathbb{A}^d_q(\Lambda)$ has some kind of quasi-commutativity, and the associativity   comes from combinatorial properties of certain polynomials appeared in the quantized cohomological Hall algebra $\mathcal{H}^d_q$ of the $d$-loop quiver. It turns out that $\mathbb{A}^d_q(\Lambda)$ is isomorphic to  $\mathcal{H}^d_q$, thus can be viewed as a combinatorial realization for $\mathcal{H}^d_q$.
\end{abstract}

\section{Introduction}

Cohomological Hall algebra of a quiver (with potential) was introduced by Kontesivech and Soibelman in \cite{KS}, which became a central tool in the study of quantized Donaldson-Thomas invariants. The underlying space of the cohomological Hall algebra of a finite quiver  is the equivariant cohomology group of representations of the quiver, while multiplication is  defined by using pull-back and push-forward operations. The discovery of cohomological Hall algebras opens new directions for motivic Donaldson-Thomas invariants, mathematical physics, representation theory and so on. Different versions of cohomological Hall algebras have been introduced and studied, see \cite{FR, SS, SV, YZ} for instance. To find a simple algebraic description of a cohomological Hall algebra is a hot topic (cf. \cite{FR19,SV}).

Let $d$ be a positive integer and $Q$ the $d$-loop quiver. The cohomological Hall algebra $\mathcal{H}^d$ of $Q$ admits an algebraic description. Namely, it is isomorphic to the space of symmetric polynomials in all possible numbers of variables endowed with a shuffle-type algebra structure. In order to develop an explicit and combinatorial setup for the study of the quantized Donaldson-Thomas invariants of $Q$,  Reineke \cite{R} introduced the quantized cohomological Hall algebra $\mathcal{H}^d_q$ and the degenerate cohomological Hall algebra $\mathcal{H}^d_0$ of $Q$. In particular, the cohomological Hall algebra $\mathcal{H}^d$ is the specialization of $\mathcal{H}^d_q$ at $q=1$, while the degenerate cohomological Hall algebra $\mathcal{H}^d_0$ is the specialization of $\mathcal{H}^d_q$ at $q=0$. Surprisingly, the degenerate cohomological Hall algebra $\mathcal{H}_0^d$ admits a combinatorial realization via partitions, which plays a key role in the study of quantized Donaldson-Thomas invariants for $Q$ (cf. \cite{R}).

The aim of this note is to pursue a similar combinatorial realization of $\mathcal{H}^d_q$ as $\mathcal{H}^d_0$ and to understand the associativity of the multiplication of $\mathcal{H}^d_q$. To do this, we introduce an algebra $\mathbb{A}^d_q(\Lambda)$ associated to partitions, which can be viewed as a quantization of the algebra $\mathbb{A}^d(\Lambda)$ introduced in \cite[Definition 5.2]{R}.  The associativity of the multiplication of $\mathbb{A}^d_q(\Lambda)$ comes from combinatorial properties of a class of polynomials, which have already appeared in the algebraic description of (quantized) cohomological Hall algebra $\mathcal{H}^d_q$.
We show that
$\mathbb{A}^d_q(\Lambda)$ is isomorphic to the quantized cohomological Hall algebra $\mathcal{H}^d_q$ and the isomorphism specializes to the isomorphism between $\mathbb{A}^d(\Lambda)$ and the degenerate cohomological Hall algebra $\mathcal{H}^d_0$.

Throughout this note, we fix a positive integer $d$.

\section{Quantized cohomological Hall algebras}

In this section, we recall the algebraic descriptions of the cohomological Hall algebra $\mathcal{H}^d$, the quantized cohomological Hall algebra $\mathcal{H}^d_q$ and the degenerate cohomological Hall algebra $\mathcal{H}^d_0$ of the $d$-loop quiver.

\subsection{$\mathcal{H}^d$ of the $d$-loop quiver}

For a non-negative integer $n$, denote by $S_n$ the symmetric group on $n$ letters. Let $\mathbb{Q}[x_1,\ldots, x_n]$ be the ring of polynomials in variables $x_1, \ldots, x_n$ with rational coefficients. Denote by $\mathbb{Q}[x_1,\ldots, x_n]^{S_n}$ the subspace of $\mathbb{Q}[x_1,\ldots, x_n]$ consisting of symmetric polynomials.

 Let
\[\mathcal{H}^d=\bigoplus_{n\geq 0}\mathbb{Q}[x_1, \ldots, x_n]^{S_n}.
\]
We endow $\mathcal{H}^d$ with a multiplication $\ast$ as follows.
Given $f_1\in \mathbb{Q}[x_1,\ldots, x_{n_1}]^{S_{n_1}}$, $f_2\in \mathbb{Q}[x_1,\ldots, x_{n_2}]^{S_{n_2}}$, let
\begin{align}
&f_1\ast f_2(x_1,\ldots,x_{n_1+n_2})\nonumber\\
&=\sum_{\{i_1,\ldots,i_{n_1}\}} f_1(x_{i_1},\ldots,x_{i_{n_1}})f_2(x_{j_1},\ldots,x_{j_{n_2}})
\prod_{l=1}^{n_2}\prod_{t=1}^{n_1}
(x_{j_l}-x_{i_t})^{d-1},\label{qe0}
\end{align}
where     the sum ranges over all the subsets $\{i_1,\ldots,i_{n_1}\}$    of $\{1,2,\ldots,n_1+n_2\}$ and $\{j_1,\ldots,j_{n_2}\}=\{1,2,\ldots,n_1+n_2\}\setminus
\{i_1,\ldots,i_{n_1}\}$. According to \cite[Theorem 2.2]{KS}, the \textit{cohomological Hall algebra} of the $d$-loop quiver is isomorphic to $(\mathcal{H}^d, \ast)$.

\subsection{$\mathcal{H}^d_q$ and $\mathcal{H}^d_0$ of $d$-loop quiver}

Let $q$ be an indeterminate and $\mathbb{Q}[q]$ the ring of polynomials in $q$ with rational coefficients.
The (naive) {\it quantized cohomological Hall algebra} $\mathcal{H}_q^d$ is the graded $\mathbb{Q}[q]$-module
\[\bigoplus_{n\geq 0}\mathbb{Q}[q][x_1,\ldots,x_n]^{S_n}
\]
endowed with the multiplication
\begin{align}
&f_1\ast f_2(x_1,\ldots,x_{n_1+n_2})\nonumber\\
&=\sum_{\{i_1,\ldots,i_{n_1}\}} f_1(x_{i_1},\ldots,x_{i_{n_1}})f_2(x_{j_1},\ldots,x_{j_{n_2}})
\prod_{l=1}^{n_2}\prod_{t=1}^{n_1}
(x_{j_l}-qx_{i_t})^{d-1},\label{qe1}
\end{align}
where $f_1\in \mathbb{Q}[q][x_1,\ldots, x_{n_1}]^{S_{n_1}}$, $f_2\in \mathbb{Q}[q][x_1,\ldots, x_{n_2}]^{S_{n_2}}$ and  the sum ranges over all the subsets $\{i_1,\ldots,i_{n_1}\}$    of $\{1,2,\ldots,n_1+n_2\}$ and $\{j_1,\ldots,j_{n_2}\}=\{1,2,\ldots,n_1+n_2\}\setminus
\{i_1,\ldots,i_{n_1}\}$, see \cite[Definition 5.3]{R}.

We can specialize the algebra $\mathcal{H}^d_q$ to any $q\in \mathbb{Q}$.  The {\it degenerate cohomological Hall algebra} $\mathcal{H}^d_0$ of the $d$-loop quiver is the specialization $\mathcal{H}^d_q|_{q=0}$ of $\mathcal{H}^d_q$ at $q=0$. In particular, $\mathcal{H}^d_0$ has the same underlying space as $\mathcal{H}^d$ and its multiplication  reduces to
\begin{align*}
&f_1\ast f_2(x_1,\ldots,x_{n_1+n_2})\nonumber\\
&=\sum_{\{i_1,\ldots,i_{n_1}\}} f_1(x_{i_1},\ldots,x_{i_{n_1}})f_2(x_{j_1},\ldots,x_{j_{n_2}})
\left(\prod_{l=1}^{n_2}
x_{j_l}\right)^{(d-1)n_1}.
\end{align*}

\subsection{A combinatorial realization of $\mathcal{H}^d_0$}

Let $\mathbb{N}$ be the set of non-negative integers.
Let $n$ be a non-negative integer. A partition of length $n$ is a vector $\lambda=(\lambda_1,\lambda_2,\ldots, \lambda_n)\in \mathbb{N}^n$ such that $\lambda_1\leq \lambda_2\leq \cdots\leq \lambda_n$.  For a partition $\lambda$, we denote by $|\lambda|$ the length of $\lambda$. For an integer vector $\mu\in \mathbb{N}^n$, denote by $\mu_{\leq}$ the partition obtained by rearranging the sequence of $\mu$ increasingly.

Let $\Lambda_n$ be the set of partitions  of length $n$. Denote by $\Lambda$ the disjoint union of all $\Lambda_n$ for $n\geq 0$. Reineke \cite[Definition 5.2]{R} introduced an associative algebra $\mathbb{A}^d(\Lambda)$, which is the $\mathbb{Q}$-vector space with basis elements $\lambda\in \Lambda$ endowed with a multiplication as follows
\begin{align}
\mu\ast \nu=(\mu_1,\ldots, \mu_m, \nu_1+m(d-1), \ldots, \nu_n+m(d-1))_{\leq},
\end{align}
where $|\mu|=m$ and $|\nu|=n$. The multiplication of $\mathbb{A}^d(\Lambda)$ is obviously associative. It has been proved by \cite[Proposition 5.4]{R} that $\mathbb{A}^d(\Lambda)$ is isomorphic to $\mathcal{H}^d_0$.

\section{The  algebra $\mathbb{A}^d_q(\Lambda)$}

In this section, we will introduce a  $\mathbb{Q}[q]$-algebra associated to partitions, which can be viewed as a quantization of $\mathbb{A}^d(\Lambda)$.
\subsection{The polynomial $g_{m,n}$}
Let $x_1,x_2,\ldots$ be infinitely many formal variables. For a non-negative integer vector ${\bf a}=(a_1,\ldots, a_n)\in \mathbb{N}^n$, we write
\[x^{\bf a}:=\prod_{i=1}^nx_{i}^{a_i}.
\]
For a pair of positive integers $(m,n)$, denote by  $g_{m,n}\in \mathbb{Q}[q][x_1,\ldots, x_{m+n}]$ the following polynomial
\begin{eqnarray}\label{gmn}
g_{m,n}:=
\prod_{s=m+1}^{m+n}\prod_{t=1}^m(x_s-qx_t)^{d-1}
=\sum_{\tiny \begin{array}{l}{\bf a}\in \mathbb{N}^m \\
{\bf b}\in \mathbb{N}^n\end{array}}c_{{m,n}}^{{\bf a},
{\bf b}}x^{({\bf a}, {\bf b})},
\end{eqnarray}
where $c_{m,n}^{{\bf a}, {\bf b}}\in \mathbb{Q}[q]$ is the coefficient of the monomial $x^{({\bf a},{\bf b})}:=x^{(a_1,\ldots,a_m,b_1,\ldots,b_n)}$ in the expansion of $g_{m,n}$. Let $l$ be a positive integer, we also denote 
\[g_{m,n}[l]:=\prod_{s=l+m+1}^{l+m+n}\prod_{t=l+1}^{l+m}(x_s-qx_t)^{d-1}.
\]

By definition, it is easy to see that $g_{m,n}$ is symmetric with respect to
$x_1,\ldots, x_m$ and $x_{m+1}, \ldots, x_{m+n}$, respectively.
Consequently, for any ${\bf a}\in \mathbb{N}^m$,
${\bf b}\in \mathbb{N}^n$
and   any permutations $\tau\in S_m,\sigma\in S_n$,
\[c_{m,n}^{\tau({\bf a}), \sigma({\bf b})}=c_{m,n}^{{\bf a}, {\bf b}}.\]

\begin{lemma}\label{l:cofficient-equality}
Let $l,m,n$ be positive integers and ${\bf a}\in \mathbb{N}^l, {\bf b}\in \mathbb{N}^m,{\bf c}\in \mathbb{N}^n$. We have
\begin{eqnarray}\label{f:coefficient-equality}
\sum_{\tiny \begin{array}{l}{\bf a}_1+{\bf a}_2={\bf a}\\ {\bf b}_1+{\bf b}_2={\bf b}
\end{array}}c_{l,m}^{{\bf a}_1, {\bf b}_1}c_{l+m,n}^{({\bf a}_2,{\bf b}_2), {\bf c}}
=\sum_{\tiny \begin{array}{l}{\bf c}_1+{\bf c}_2={\bf c}\\ {\bf b}_1+{\bf b}_2={\bf b}
\end{array}}c_{l,m+n}^{{\bf a}, ({\bf b}_2, {\bf c}_2)}c_{m,n}^{{\bf b}_1, {\bf c}_1},
\end{eqnarray}
where the sums are taking over all non-negative integer vectors ${\bf a}_1$, ${\bf a}_2$, ${\bf b}_1$, ${\bf b}_2$, ${\bf c}_1$, ${\bf c}_2$.
\end{lemma}

\begin{proof}
By definition, we have
\begin{eqnarray*}
g_{l,m}g_{l+m,n}&=&\prod_{s=l+1}^{l+m}\prod_{t=1}^l(x_s-qx_t)^{d-1} \prod_{s=l+m+1}^{l+m+n}\prod_{t=1}^{l+m}(x_s-qx_t)^{d-1} \\[5pt]
&=&\left(\sum_{\tiny \begin{array}{l}{\bf a}_1\in\mathbb{N}^l \\
{\bf b}_1\in\mathbb{N}^m\end{array}}c_{{l,m}}^{{\bf a}_1,{\bf b}_1}
x^{({\bf a}_1, {\bf b}_1)}\right) \left(\sum_{\tiny \begin{array}{l}
{\bf a}_2\in \mathbb{N}^l \\ {\bf b}_2\in \mathbb{N}^m\\ {\bf c}
\in \mathbb{N}^n\end{array}}c_{{l+m,n}}^{({\bf a}_2,{\bf b}_2),
{\bf c}}x^{({\bf a}_2, {\bf b}_2, {\bf c})}\right)\\[5pt]
&=&\sum_{\tiny \begin{array}{l}{\bf a}_1,{\bf a}_2\in \mathbb{N}^l \\ {\bf b}_1, {\bf b}_2\in \mathbb{N}^m\\ {\bf c}\in \mathbb{N}^n\end{array}}c_{{l,m}}^{{\bf a}_1,{\bf b}_1}c_{{l+m,n}}^{({\bf a}_2,{\bf b}_2),{\bf c}}x^{({\bf a}_1+{\bf a}_2, {\bf b}_1+{\bf b}_2, {\bf c})}.
\end{eqnarray*}
Thus the left hand side of (\ref{f:coefficient-equality}) is  the coefficient of the monomial $x^{({\bf a}, {\bf b}, {\bf c})}$ in the expansion of $g_{l,m}g_{l+m,n}$.

On the other hand, consider the  following polynomial
\begin{eqnarray*}
g_{l,m+n}(g_{m,n}[l])&=&\prod_{s=l+1}^{l+m+n}\prod_{t=1}^l (x_s-qx_t)^{d-1} \prod_{s=l+m+1}^{l+m+n}\prod_{t=l+1}^{l+m}(x_s-qx_t)^{d-1}\\[5pt]
&=&\left(\sum_{\tiny \begin{array}{l}{\bf a}\in \mathbb{N}^l\\ {\bf b}_1\in \mathbb{N}^m\\ {\bf c}_1\in \mathbb{N}n\end{array}}c_{l,m+n}^{{\bf a}, ({\bf b}_1,{\bf c}_1)}x^{({\bf a}, {\bf b}_1,{\bf c}_1)}\right)\left(\sum_{\tiny \begin{array}{l}{\bf b}_2\in \mathbb{N}^m\\ {\bf c}_2\in \mathbb{N}n\end{array}}c_{m,n}^{{\bf b}_2,{\bf c}_2}x^{({\bf 0}, {\bf b}_2,{\bf c}_2)}\right)\\[5pt]
&=&\sum_{\tiny \begin{array}{l}{\bf a}\in \mathbb{N}^l\\ {\bf b}_1,{\bf b}_2\in \mathbb{N}^m\\ {\bf c}_1,{\bf c}_2\in \mathbb{N}n\end{array}}c_{l,m+n}^{{\bf a}, ({\bf b}_1,{\bf c}_1)}c_{m,n}^{{\bf b}_2,{\bf c}_2}x^{({\bf a}, {\bf b}_1+{\bf b}_2,{\bf c}_1+{\bf c}_2)}.
\end{eqnarray*}
Then the right hand side of (\ref{f:coefficient-equality}) is  the coefficient of the monomial $x^{({\bf a}, {\bf b}, {\bf c})}$ in the expansion of $g_{l,m+n}(g_{m,n}[l])$. The equality (\ref{f:coefficient-equality}) follows from the following obvious relation
\[g_{l,m}g_{l+m,n}=g_{l,m+n}(g_{m,n}[l]).\]
\end{proof}

\subsection{An auxiliary algebra $\mathbb{A}^d_q(\hat{\Lambda})$}
Denote by $\hat{\Lambda}$ the set of disjoint union of all $\mathbb{N}^n$ for $n\geq 0$. Let $\mathbb{A}^d_q(\hat{\Lambda})$ be the free $\mathbb{Q}[q]$-module with basis elements $\lambda\in \hat{\Lambda}$. For an element $\lambda\in \hat{\Lambda}$, we denote by $|\lambda|$ the length of $\lambda$, that is, $\lambda\in \mathbb{N}^{|\lambda|}$. Denote by $\emptyset$ the unique element of $\mathbb{N}^0$.
Let $\mu\in \mathbb{N}^m, \nu\in \mathbb{N}^n$ with $m,n\geq 1$, we define
$
\mu\ast \emptyset=\mu=\emptyset \ast \mu
$ and
\begin{align}\label{auxi}
\mu\ast\nu=\sum_{\tiny\begin{array}{l}{\bf a}\in \mathbb{N}^m\\ {\bf b}\in \mathbb{N}^n \end{array}}c_{m,n}^{{\bf a}, {\bf b}}(\mu+{\bf a}, \nu+{\bf b}),
\end{align}
where $c_{m,n}^{{\bf a}, {\bf b}}$ is as defined in \eqref{gmn}. Extend the multiplication defined in \eqref{auxi} $\mathbb{Q}[q]$-linearly to a multiplication of $\mathbb{A}^d_q(\hat{\Lambda})$.

\begin{lemma}
$(\mathbb{A}^d_q(\hat{\Lambda}), \ast)$ is an associative algebra with unit.
\end{lemma}
\begin{proof}
By definition $\emptyset$ is the unit. It remains to prove the associativity, that is, for $\mu,\nu,w\in  \hat{\Lambda}$, one has  $(\mu\ast \nu)\ast w=\mu\ast(\nu\ast w)$. There is nothing to prove if one of $\mu,\nu,w$ is $\emptyset$. Without loss of generality, we may assume that $|\mu|=l, |\nu|=m,|w|=n \geq 1$. By definition,
\begin{align*}
(\mu\ast \nu)\ast w&=\sum_{\tiny \begin{array}{l} {\bf a}_1\in \mathbb{N}^l\\ {\bf b}_1\in \mathbb{N}^m\end{array}}c_{l,m}^{{\bf a}_1, {\bf b}_1}(\mu+{\bf a}_1, \nu+{\bf b}_1)\ast w\\
&=\sum_{\tiny \begin{array}{l} {\bf a}_1\in \mathbb{N}^l\\ {\bf b}_1\in \mathbb{N}^m\end{array}}c_{l,m}^{{\bf a}_1, {\bf b}_1} \sum_{\tiny \begin{array}{l} {\bf a}_2\in \mathbb{N}^l\\ {\bf b}_2\in \mathbb{N}^m\\ {\bf c}\in \mathbb{N}^n\end{array}}c_{l+m,n}^{({\bf a}_2, {\bf b}_2), {\bf c}}(\mu+{\bf a}_1+{\bf a_2}, \nu+{\bf b}_1+{\bf b}_2, w+{\bf c})\\
&=\sum_{\tiny \begin{array}{l} {\bf a}_1, {\bf a}_2\in \mathbb{N}^l\\ {\bf b}_1,{\bf b}_2\in \mathbb{N}^m\\ {\bf c}\in \mathbb{N}^n\end{array}}c_{l,m}^{{\bf a}_1, {\bf b}_1}
c_{l+m,n}^{({\bf a}_2, {\bf b}_2), {\bf c}}
(\mu+{\bf a}_1+{\bf a_2}, \nu+{\bf b}_1+{\bf b}_2, w+{\bf c}).
\end{align*}
Therefore, the coefficient of $(\mu+{\bf a}, \nu+{\bf b}, w+{\bf c})$ in $(\mu\ast\nu)\ast w$ is
\[\sum\limits_{\tiny \begin{array}{l}{\bf a}_1+{\bf a}_2={\bf a}\\ {\bf b}_1+{\bf b}_2={\bf b}
\end{array}}c_{l,m}^{{\bf a}_1, {\bf b}_1}c_{l+m,n}^{({\bf a}_2,{\bf b}_2), {\bf c}}.\]

On the other hand, we have
\begin{eqnarray*}
\mu\ast(\nu\ast w)&=&\mu\ast \sum_{\tiny \begin{array}{l}{\bf b}_2\in \mathbb{N}^n\\ {\bf c}_2\in \mathbb{N}^n\end{array}}c_{m,n}^{{\bf b}_2, {\bf c}_2}(\nu+{\bf b}_2, w+{\bf c}_2)\\
&=&\sum_{\tiny \begin{array}{l} {\bf a}\in \mathbb{N}^l\\ {\bf b}_1,{\bf b}_2\in \mathbb{N}^m\\ {\bf c}_1, {\bf c}_2\in \mathbb{N}^n\end{array}}c_{l,m+n}^{{\bf a}, ({\bf b}_2, {\bf c}_2)}
c_{m,n}^{{\bf b}_1, {\bf c}_1}
(\mu+{\bf a}, \nu+{\bf b}_1+{\bf b}_2, w+{\bf c}_1+{\bf c}_2).
\end{eqnarray*}
Thus the coefficient of $(\mu+{\bf a}, \nu+{\bf b}, w+{\bf c})$ in $\mu\ast(\nu\ast w)$ is
\[\sum\limits_{\tiny \begin{array}{l}{\bf c}_1+{\bf c}_2={\bf c}\\ {\bf b}_1+{\bf b}_2={\bf b}
\end{array}}c_{l,m+n}^{{\bf a}, ({\bf b}_2, {\bf c}_2)}c_{m,n}^{{\bf b}_1, {\bf c}_1}.\]
We conclude that $(\mu\ast\nu)\ast w=\mu\ast(\mu\ast w)$ by Lemma \ref{l:cofficient-equality}.
\end{proof}
\begin{remark}
In fact,  $\mathbb{A}^d_q(\hat{\Lambda})$ is an $\mathbb{N}$-graded algebra
\[\mathbb{A}^d_q(\hat{\Lambda})=\bigoplus_{n\in \mathbb{N}}\mathbb{A}^d_q(\hat{\Lambda})_n,
\]
whose $n$-th component is the free $\mathbb{Q}[q]$-module with basis $\mathbb{N}^n$.
\end{remark}

\subsection{The algebra $\mathbb{A}^d_q(\Lambda)$ associated to partitions}

Let $\mu,\nu\in \hat{\Lambda}$, we denote by $\mu\sim \nu$ if $\mu$ is a permutation of $\nu$. Clearly, $\sim$ is an equivalence relation of $\hat{\Lambda}$ and
we may identify $\Lambda$ with the quotient $\hat{\Lambda}/\sim$
 of $\hat{\Lambda}$ by $\sim$. Let $W$ be the
  $\mathbb{Q}[q]$-subspace of $\mathbb{A}^d_q(\hat{\Lambda})$ spanned by
\[
\{\mu\circleddash\nu~|~ \text{for all $\mu,\nu\in \hat{\Lambda}$
such that $\mu\sim \nu$}\}.
\]
Here, to
distinguish with the minus of vectors, we use $\circleddash$ to represent the minus
of $\mathbb{A}^d_q(\hat{\Lambda})$.

\begin{lemma}\label{l:homogeneous}
The two-side ideal $( W)$ of
 $\mathbb{A}^d_q(\hat{\Lambda})$ generated by $W$ is exactly $W$. 
\end{lemma}
\begin{proof}
Let $\mu\in \hat{\Lambda}$ with $|\mu|=m$ and $\sigma\in S_m$. For any $w\in \hat{\Lambda}$ with $|w|=n$, we have
\begin{align*}%\label{e:ideal}
&(\mu\circleddash\sigma(\mu))\ast w\\
&=
\sum_{\tiny\begin{array}{l}{\bf a}\in\mathbb{N}^m\\
{\bf b}\in\mathbb{N}^n\end{array}}c_{m,n}^{{\bf a}, {\bf b}}
(\mu+{\bf a}, w+{\bf b})\circleddash\sum_{\tiny \begin{array}{l}{\bf a}\in \mathbb{N}^m\\
{\bf b}\in \mathbb{N}^n\end{array}}c_{m,n}^{{\bf a}, {\bf b}}
(\sigma(\mu)+{\bf a}, w+{\bf b})\\
&=\sum_{\tiny\begin{array}{l}{\bf a}\in\mathbb{N}^m\\
{\bf b}\in\mathbb{N}^n\end{array}}c_{m,n}^{{\bf a}, {\bf b}}
(\mu+{\bf a}, w+{\bf b})\circleddash\sum_{\tiny \begin{array}{l}{\bf a}\in \mathbb{N}^m\\
{\bf b}\in \mathbb{N}^n\end{array}}c_{m,n}^{{\bf a}, {\bf b}}
(\sigma(\mu)+\sigma({\bf a}), w+{\bf b})\\
&=\sum_{\tiny\begin{array}{l}{\bf a}\in\mathbb{N}^m\\
{\bf b}\in\mathbb{N}^n\end{array}}c_{m,n}^{{\bf a}, {\bf b}}[(\mu+{\bf a},w+{\bf b})\circleddash (\sigma(\mu+{\bf a}), w+{\bf b})].
\end{align*} 
As a consequence, $(\mu\circleddash\sigma(\mu))\ast w$ belongs to $W$. Similarly, one can show that $w\ast (\mu\circleddash\sigma(\mu))$ belongs to $W$ as well.
\end{proof}

\begin{definition}\label{d:quantized-alg}
The algebra $\mathbb{A}^d_q(\Lambda)$ associated to partitions is the quotient algebra $\mathbb{A}^d_q(\hat{\Lambda})/(W)$.
\end{definition}

The following result is a direct consequence of Lemma \ref{l:homogeneous}, Definition \ref{d:quantized-alg} and (\ref{auxi}).

\begin{proposition}
\begin{enumerate}
\item[(1)] $\mathbb{A}^d_q(\Lambda)$ is an $\mathbb{N}$-graded associative algebra with unit, whose $n$-th component  is a free $\mathbb{Q}[q]$-module with basis $\Lambda_n$;
\item[(2)]  Let $\mu\in \Lambda_m, \nu\in \Lambda_n$. The multiplication of $\mu$ with $\nu$ in  $\mathbb{A}^d_q(\Lambda)$ is as follows
\begin{eqnarray}\label{e:multiplication}
\mu\ast\nu=\sum_{\tiny \begin{array}{l}{\bf a}\in \mathbb{N}^m\\ {\bf b}\in \mathbb{N}^n\end{array}}c_{m,n}^{{\bf a}, {\bf b}}(\mu+{\bf a}, \nu+{\bf b})_{\leq},
\end{eqnarray}
where $c_{m,n}^{{\bf a}, {\bf b}}$ is as defined in \eqref{gmn}.
\end{enumerate}
\end{proposition}

Note that $c_{m,n}^{{\bf a},{\bf b}}$ is a polynomial in $q$. Hence we can specialize the value of $c_{m,n}^{{\bf a},{\bf b}}$ at $q=0$. The following result suggests that $\mathbb{A}^d_q(\Lambda)$ is a quantization of $\mathbb{A}^d(\Lambda)$.
\begin{corollary}
By specializing $q$ to $0$, we obtain
\[\mathbb{A}^d_q(\Lambda)|_{q=0}\cong \mathbb{A}^d(\Lambda).
\]
\end{corollary}
\begin{example} Assume that $d=4$.
 Let $\mu=(0,2),\nu=(1)$. Then
\[g_{2,1}=(x_3-qx_1)^3(x_3-qx_2)^3.\]
One can compute $\mu\ast\nu$ in $\mathbb{A}^4_q(\Lambda)$ as follows.
\begin{align*}
\mu\ast\nu&=(q^6+3q^4+9q^2)\cdot(1, 3, 5)+(3q^2-3q^5)\cdot(2, 2, 5)\\
&\quad+(9q^4-3q^5-10q^3)\cdot(2, 3, 4)
+3q^4\cdot(3, 3, 3)+(3q^2-q^3)\cdot(0, 4, 5)\\
&\quad -9q^3\cdot(1, 4, 4) -3q\cdot(0, 3, 6) -3q\cdot(1, 2, 6)+(0, 2, 7).
\end{align*}
\end{example}
Let $q^{-1}$ be the formal inverse of $q$ and
$\mathbb{Q}[q^{\pm}]$ the Laurent polynomials in $q$
with rational coefficients. Let us consider the tensor
 product $\mathbb{Q}[q^{\pm}]\otimes_{\mathbb{Q}[q]}
 \mathbb{A}^d_q(\Lambda)$. In particular,
 $\mathbb{A}^d_q(\Lambda)$ is a subalgebra of
  $\mathbb{Q}[q^{\pm}]\otimes_{\mathbb{Q}[q]}\mathbb{A}^d_q(\Lambda)$.
 The bar involution $\bar{ } : \mathbb{Q}[q^\pm]\to \mathbb{Q}[q^\pm]$ given by $q^i\mapsto q^{-i}$($i\in \mathbb{Z}$),  induces a bar involution of
\begin{eqnarray*}
\bar{}:&&\mathbb{Q}[q^{\pm}]\otimes_{\mathbb{Q}[q]}\mathbb{A}^d_q(\Lambda)\to \mathbb{Q}[q^{\pm}]\otimes_{\mathbb{Q}[q]}\mathbb{A}^d_q(\Lambda),\\
&&f=\sum a_\mu\mu\mapsto \bar{f}=\sum\bar{a}_\mu\mu, 
\end{eqnarray*}
where $a_\mu\in \mathbb{Q}[q^\pm]$. The following is a direct consequence of the property of the polynomial $g_{m,n}$ and (\ref{e:multiplication}).
\begin{proposition}
Let $\mu\in \Lambda_m,\nu\in \Lambda_n$. We have
\[\mu\ast\nu=q^{(d-1)mn}\overline{\nu\ast\mu}.
\]
\end{proposition}

\section{Isomorphism from $\mathbb{A}^d_q(\Lambda)$ to $\mathcal{H}^d_q$}

For a non-negative integer sequence $\lambda=(\lambda_1,\ldots, \lambda_n)\in \mathbb{N}^n$, we denote by
\begin{align}\label{eq2}
M_{\lambda}(x_1,\ldots,x_n):=\sum_{\sigma\in S_n}x_{\sigma(1)}^{\lambda_1}\cdots x_{\sigma(n)}^{\lambda_n}\in \mathcal{H}^d_q.
\end{align}
Let \[c({\lambda})=|\{\sigma\in S_n~|~ \lambda=(\lambda_{\sigma(1)}, \dots, \lambda_{\sigma(n)})\}|.
\]
Denote by $S_\lambda$ the set of different permutations of the parts of $\lambda$, that is, 
\[S_\lambda=\{(t_1,\dots, t_n)~|~(t_1,\dots, t_n)=(\lambda_{\sigma(1)},\dots, \lambda_{\sigma(n)}), \text{for some $\sigma\in S_n$}\}.
\]
Clearly, 
\[c(\lambda)=\frac{n!}{|S_\lambda|}.
\]
Let $m_\lambda(x_1,\ldots,x_n)$ be the monomial symmetric function of $\lambda$, that is, the sum of all the monomials $x_1^{\pi(1)}\cdots x_n^{\pi(n)}$, where $\pi\in S_\lambda$.
It is obvious that
\[M_{\lambda}(x_1,\ldots,x_n)=c(\lambda)\cdot m_{\lambda}(x_1,\ldots,x_n).\]
It is also easy to see that $m_{\lambda}(x_1,\ldots,x_n)=m_{\lambda'}(x_1,\ldots,x_n)$ and $M_{\lambda}(x_1,\ldots,x_n)=M_{\lambda'}(x_1,\ldots,x_n)$  for any $\lambda'\in S_{\lambda}$.

The following proposition is already  known, see, for example, Carvalho and  D'Agostino \cite{CD},  we give a short proof here for the convenience of the readers.

\begin{proposition}\label{pmp}
Let $\lambda,\mu$ be two partitions of length $n$. Then
\[M_{\lambda}(x_1,\ldots,x_n)m_{\mu}(x_1,\ldots,x_n)
=\sum_{w\in S_\mu}M_{\lambda+w}(x_1,\ldots,x_n).\]
\end{proposition}

\begin{proof} Fix $w\in S_{\mu}$, choose a permutation $\tau\in S_n$ (not necessarily unique) such that $w_i=\mu_{\tau(i)}.$ Then we have
\begin{eqnarray*}M_{\lambda+w}(x_1,\cdots,x_n)&=&\sum_{\pi\in S_{n}}x_{\pi(1)}^{\lambda_1+w_1}\cdots x_{\pi(n)}^{\lambda_n+w_n}\\
&=&\sum_{\pi\in S_{n}}x_{\pi(1)}^{\lambda_1+\mu_{\tau(1)}}\cdots x_{\pi(n)}^{\lambda_n+\mu_{\tau(n)}}.
\end{eqnarray*}
%For a given monomial $x_{\pi(1)}^{\lambda_1+\mu_{\tau(1)}}\cdots x_{\pi(n)}^{\lambda_n+\mu_{\tau(n)}}$ of $P_{\nu}(x_1,\cdots,x_n)$, we first show that this monomial will appear exactly $c(\mu)$ times in the expansion of $P_{\lambda}(x_1,\cdots,x_n)m_{\mu}(x_1,\cdots,x_n)$.
On the other hand,
\begin{align*}
M_{\lambda}(x_1,\cdots,x_n)M_{\mu}(x_1,\cdots,x_n)
&=\sum_{\pi\in S_n}x_{\pi(1)}^{\lambda_1}\cdots x_{\pi(n)}^{\lambda_n}\sum_{\beta\in S_n}x_{\beta(1)}^{\mu_1}\cdots x_{\beta(n)}^{\mu_n}\\
&=\sum_{\pi\in S_n}\left(\sum_{\beta\in S_n}x_{\pi(1)}^{\lambda_1+
\mu_{\beta^{-1}(\pi(1))}}\cdots x_{\pi(n)}^{\lambda_n+\mu_{\beta^{-1}(\pi(n))}}\right),
\end{align*}
For a given monomial $x_{\pi(1)}^{\lambda_1+\mu_{\tau(1)}}\cdots x_{\pi(n)}^{\lambda_n+\mu_{\tau(n)}}$ of $M_{\lambda+w}(x_1,\ldots,x_n)$,
there are $c(\mu)$  permutations $\beta\in S_n$  in the summand
\[\sum_{\beta\in S_n}x_{\pi(1)}^{\lambda_1+
\mu_{\beta^{-1}(\pi(1))}}\cdots x_{\pi(n)}^{\lambda_n+\mu_{\beta^{-1}(\pi(n))}}\]
such that $\mu_{\beta^{-1}(\pi(i))}=\mu_{\tau(i)}$ for $i=1,\ldots,n$. Moreover, since there are $|S_{\mu}|$ different $\lambda+w$, and each $\lambda+w$ corresponds to $c(\mu)$ different permutations $\beta\in S_n$, we see that
\begin{align}
M_{\lambda}(x_1,\ldots,x_n)M_{\mu}(x_1,\ldots,x_n)
=c(\mu)\sum_{w\in S_\mu}M_{\lambda+w}(x_1,\ldots,x_n).
\end{align}
Therefore
\[M_{\lambda}(x_1,\ldots,x_n)m_{\mu}(x_1,\ldots,x_n)
=\sum_{w\in S_\mu}M_{\lambda+w}(x_1,\ldots,x_n).\]
This completes the proof. 
\end{proof}

\begin{example}
For example, let $\lambda=(0,1,2),\mu=(1,1,3)$. Then there are 3 different permutations of $\mu$, i.e., $S_{\mu}=\{(3,1,1),(1,3,1),(1,1,3)\}$. 
Then we have
\[M_{(0,1,2)}m_{(1,1,3)}=M_{(3,2,3)}+M_{(1,4,3)}+M_{(1,2,5)}=M_{(2,3,3)}
+M_{(1,3,4)}+M_{(1,2,5)}.\]
\end{example}

\begin{corollary}\label{pm}
Let $\lambda,\mu$ be two partitions of length $n$. Then
\[m_{\lambda}(x_1,\ldots,x_n)m_{\mu}(x_1,\ldots,x_n)
=\sum_{w\in S_\mu}\frac{c(\lambda+w)}{c(\lambda)} m_{\lambda+w}(x_1,\ldots,x_n).\]
\end{corollary}

%\begin{prop}
%Let $\lambda=(0\le \lambda_1\le\cdots\le\lambda_n)$ be a partition of length $n$ and $k\ge0$ be an integer. Then
%\[
%\sum_{\{i_1,\ldots,i_n\}}P_{\lambda}(x_{i_1},\ldots,x_{i_n})
%=\frac{1}{k!}\cdot P_{\lambda'}(x_1,\ldots,x_{n+k}),
%\]
%where the sum takes over all the $n$-subsets of $\{1,2,\ldots,n+k\}$ and $\lambda'=(0,\ldots,0,\lambda_1,\ldots,\lambda_n)$ is obtained from $\lambda$ by appending $k$ zeros at the front.
%\end{prop}

\begin{proposition}\label{2pp}
Let $\lambda=(\lambda_1,\ldots,\lambda_n)\in\Lambda_n$  and $\mu=(\mu_1,\ldots,\mu_k)\in\Lambda_k$. Then
\[
\sum_{\{i_1,\ldots,i_n\}} M_{\lambda}(x_{i_1},\ldots,x_{i_n})
M_{\mu}(x_{j_1},\ldots,x_{j_k})=M_{(\lambda,\mu)_{\le}}(x_1,\ldots,x_{n+k}),
\]
where the sum takes over all the subsets $\{i_1,\ldots,i_n\}$ of $\{1,\ldots,n+k\}$ and $\{j_1,\ldots,j_k\}=\{1,\ldots,n+k\}\setminus\{i_1,\ldots,i_n\}$.
\end{proposition}

\begin{proof} By definition \eqref{eq2},
\begin{align*}
&\sum_{\{i_1,\ldots,i_n\}}M_{\lambda}(x_{i_1},\ldots,x_{i_n})
M_{\mu}(x_{j_1},\ldots,x_{j_k})\\
&=\sum_{\{i_1,\ldots,i_n\}}\left(\sum_{\sigma\in S_n}x_{\sigma(i_1)}^{\lambda_1}\cdots x_{\sigma(i_n)}^{\lambda_n}\right)\left(\sum_{\pi\in S_k}x_{\pi(j_1)}^{\mu_1}\cdots x_{\pi(j_k)}^{\mu_k}\right)\\
&=\sum_{(\{i_1,\ldots,i_n\},\sigma,\pi)}x_{\sigma(i_1)}^{\lambda_1}\cdots x_{\sigma(i_n)}^{\lambda_n}x_{\pi(j_1)}^{\mu_1}\cdots x_{\pi(j_k)}^{\mu_k}.
\end{align*}
It is easy to see that the set of triples $\{(\{i_1,\ldots,i_n\},\sigma,\pi)\}$
is in bijection with the set of permutations $S_{n+k}$. More precisely, a triple $(\{i_1,\ldots,i_n\},\sigma,\pi)$ can be associated to the permutation \[\tau=(\tau(1),\ldots,\tau(n+k))
:=(\sigma(i_1),\ldots,\sigma(i_n),\pi(j_1),\ldots,\pi(j_k)).\]
 Thus we have
\begin{align*}
&\sum_{\{i_1,\ldots,i_n\}}M_{\lambda}(x_{i_1},\ldots,x_{i_n})
M_{\mu}(x_{j_1},\ldots,x_{j_k})\\
&=\sum_{\tau\in S_{n+k}}x_{\tau(1)}^{\lambda_1}\cdots x_{\tau(n)}^{\lambda_n}x_{\tau(n+1)}^{\mu_1}\cdots x_{\tau(n+k)}^{\mu_k}\\
&=M_{(\lambda,\mu)}
(x_1,\ldots,x_{n+k})\\
&=M_{(\lambda,\mu)_{\leq}}(x_1,\ldots,x_{n+k}).
\end{align*}
This completes the proof.
\end{proof}

\begin{theorem}\label{t:main-thm}
Let $\Phi:\mathbb{A}^d_q(\Lambda)\to \mathcal{H}^d_q$ be the $\mathbb{Q}[q]$-linear map determined by
\begin{align}\label{maph}
\Phi(\lambda)=M_\lambda(x_1,\ldots, x_n),
\end{align}
where $\lambda$ is a partition of length $n$.
Then $\Phi$ is an algebra isomorphism between $\mathbb{A}^d_q(\Lambda)$ and $\mathcal{H}^d_q$.
\end{theorem}

\begin{proof}
Since $M_{\lambda}=c(\lambda)\cdot m_{\lambda}$,  $\{M_{\lambda}|\lambda\in\Lambda\}$ forms a basis of $\mathcal{H}_q^d$. It is easy to see that $\Phi$ is a vector space isomorphism. To show that $\Phi$ is an algebra isomorphism, we need to show $\Phi(\mu\ast\nu)=\Phi(\mu)\ast \Phi(\nu)$, that is,
\begin{eqnarray}\label{pp}
M_{\mu}\ast M_{\nu}=\sum_{\tiny\begin{array}{l}{\bf a}\in \mathbb{N}^n\\ {\bf b}\in \mathbb{N}^k\end{array}}c_{n,k}^{{\bf a}, {\bf b}}\cdot M_{(\mu+{\bf a}, \nu+{\bf b})_{\leq }},
\end{eqnarray}
for any two partitions $\mu=(\mu_1,\ldots,\mu_n)$ and $\nu=(\nu_1,\ldots,\nu_k)$.

By definition  \eqref{qe1},
\begin{align*}
&M_{\mu}\ast M_{\nu}(x_1,\ldots,x_{n+k})\\
&=\sum_{\{i_1,\ldots,i_n\}} M_{\mu}(x_{i_1},\ldots,x_{i_{n}})M_{\nu}(x_{j_1},\ldots,x_{j_{k}})
\prod_{l=1}^{k}\prod_{t=1}^{n}
(x_{j_l}-qx_{i_t})^{d-1}.
\end{align*}
Denote by $I=\{i_1,\ldots, i_n\}\subset \{1,\ldots, n+k\}$ and $J=\{j_1,\ldots, j_k\}$ such that $I\cup J=\{1,\ldots, n+k\}$.
Let \[g_{n,k}(I,J):=\prod_{s=1}^k\prod_{t=1}^n(x_{j_s}-qx_{i_t})^{d-1}.\] Clearly, $g_{n,k}(I,J)$ is symmetric with respect to $x_{i_1}, \ldots, x_{i_n}$ and $x_{j_1}, \ldots, x_{j_k}$, respectively. Hence we have
\begin{align*}
g_{n,k}(I,J)&=\sum_{\tiny\begin{array}{l}{\bf a}\in \mathbb{N}^n\\ {\bf b}\in \mathbb{N}^k\end{array}}c_{n,k}^{{\bf a}, {\bf b}}x_{i_1}^{a_1}\cdots x_{i_n}^{a_n}x_{j_1}^{b_1}\cdots x_{j_k}^{b_k}\\
&=\sum_{\tiny\begin{array}{l}{\bf a}\in \Lambda_n\\ {\bf b}\in \Lambda_k\end{array}}c_{n,k}^{{\bf a}, {\bf b}}m_{\bf a}(x_{i_1}, \ldots, x_{i_n})m_{\bf b}(x_{j_1},\ldots, x_{j_k}).
\end{align*}
Therefore,
\begin{align*}
&M_{\mu}(x_{i_1},\ldots, x_{i_n})M_\nu(x_{j_1},\ldots, x_{j_k})g_{n,k}(I,J)\\[5pt]
&=\sum_{\tiny\begin{array}{l}{\bf a}\in \Lambda_n\\ {\bf b}\in \Lambda_k\end{array}}c_{n,k}^{{\bf a}, {\bf b}}M_\mu(x_{i_1}, \ldots, x_{i_n})m_{\bf a}(x_{i_1}, \ldots, x_{i_n})M_\nu(x_{j_1}, \ldots, x_{j_k})m_{\bf b}(x_{j_1},\ldots, x_{j_k})\\[5pt]
&=\sum_{\tiny\begin{array}{l}{\bf a}\in \Lambda_n\\ {\bf b}\in \Lambda_k\end{array}}c_{n,k}^{{\bf a}, {\bf b}}\sum_{\lambda\in S_{\bf a}}M_{(\mu+\lambda)_\leq}(x_{i_1}, \ldots, x_{i_n})\sum_{w\in S_{\bf b}}M_{(\nu+w)_\leq}(x_{j_1}, \ldots, x_{j_k}) \\[5pt]
&=\sum_{\tiny\begin{array}{l}{\bf a}\in \mathbb{N}^n\\ {\bf b}\in \mathbb{N}^k\end{array}}c_{n,k}^{{\bf a}, {\bf b}}M_{(\mu+{\bf a})_\leq}(x_{i_1}, \ldots, x_{i_n})M_{(\nu+{\bf b})_\leq}(x_{j_1}, \ldots, x_{j_k}),
\end{align*}
where the second equality follows from  Proposition \ref{pmp}.
Consequently,
\begin{align*}
&M_\mu\ast M_\nu(x_1,\ldots, x_{n+k})\\
&=\sum_{\{i_1,\ldots,i_n\}} M_{\mu}(x_{i_1},\ldots,x_{i_{n}})M_{\nu}(x_{j_1},\ldots,x_{j_{k}})
g_{n,k}(I,J)\\[5pt]
&=\sum_{\tiny\begin{array}{l}{\bf a}\in \mathbb{N}^n\\ {\bf b}\in \mathbb{N}^k\end{array}}c_{n,k}^{{\bf a}, {\bf b}}\sum_{\{i_1,\ldots, i_n\}}M_{(\mu+{\bf a})_\leq}(x_{i_1},\ldots, x_{i_n})M_{(\nu+{\bf b})_\leq}(x_{j_1}, \ldots, x_{j_k})\\[5pt]
&=\sum_{\tiny\begin{array}{l}{\bf a}\in \mathbb{N}^n\\ {\bf b}\in \mathbb{N}^k\end{array}}c_{n,k}^{{\bf a}, {\bf b}}M_{(\mu+{\bf a},\nu+{\bf b})_{\leq}},
\end{align*}
where the last equality follows from Proposition \ref{2pp}.
\end{proof}
To conclude, we remark that in the special case $d=2$,  the product (\ref{e:multiplication}) can be defined combinatorially.
Given  $\mu\in \Lambda_n$ and $\nu\in \Lambda_k$, let $G_{\mu,\nu}$ denote the set of  direct bipartite graphs between nodes set $U=\{u_1,\ldots,u_n\}$ and $V=\{v_1,\ldots,v_k\}$ with an orientation of each edge.  Clearly, $G_{\mu,\nu}$ has $2^{n\cdot k}$ such direct bipartite graphs.
Given a    bipartite graph $G\in G_{\mu,\nu}$, %let
%\[wt(e)=-q(m-1),\]
%if $e=v_j\rightarrow u_i$ is an edge pointed from  $v_j$  to $u_i$, otherwise  let \[wt(e)=m-1.\]
for $1\le i\le n$ and $1\le j\le k$,
let
\begin{align*}
a_i&=\mu_i+ \#\{\text{directed edges pointed to the node $u_i$}\},\\ 
b_j&=\nu_j+ \#\{\text{directed edges pointed to the node $v_j$}\}.
\end{align*}
Construct a vector $h(G)$ of $\mathbb{A}^2_q(\Lambda)$ according to $G$ as follows. The partition of $h(G)$ is  $(a_1,\ldots,a_n,b_1,\ldots,b_k)_{\le}$. The coefficient of $h(G)$ is
$(-q)^{m(G)}$,
where
\[m(G)=\sum_{i=1}^n\#\{\text{directed edges pointed to the node $u_i$}\}.\]
Now we can define
\begin{align*}\label{mp}
\mu\ast\nu=\sum_{G\in G_{\mu,\nu}}h(G).
\end{align*}

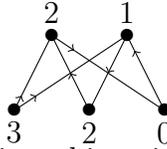
\begin{figure}[ht]
\begin{center}
\begin{tikzpicture}
\node at (0mm,0mm){$\bullet$};  \node [below] (0mm,0mm) {3};
\node at (10mm,0mm){$\bullet$}; \node at (10mm,-3mm) {2};
\node at (20mm,0mm){$\bullet$}; \node at (20mm,-3mm) {0};
\node at (5mm,10mm){$\bullet$}; \node at (5mm,13mm) {2};
\node at (15mm,10mm){$\bullet$};\node at (15mm,13mm) {1};

\draw[->](0mm,0mm)--(1mm,2mm);\draw(1mm,2mm)--(5mm,10mm);
\draw[->](0mm,0mm)--(3mm,2mm);\draw(3mm,2mm)--(15mm,10mm);

\draw[->](10mm,0mm)--(7.5mm,5mm);\draw(7.5mm,5mm)--(5mm,10mm);
\draw[->](15mm,10mm)--(12.5mm,5mm);\draw(12.5mm,5mm)--(10mm,0mm);

\draw[->](5mm,10mm)--(8mm,8mm);\draw(8mm,8mm)--(20mm,0mm);
\draw[->](20mm,0mm)--(16mm,8mm);\draw(16mm,8mm)--(15mm,10mm);

\end{tikzpicture}
\end{center}
\vspace{-.8cm}
\caption{A direct bipartite graph $G$ in $G_{\mu,\nu}$.}\label{PPP}
\end{figure}

For example, let $\mu=(1,2)$ and $\nu=(0,2,3)$. Then $G_{\mu,\nu}$ has $2^6$ direct bipartite graphs, we illustrate one of them as in Figure \ref{PPP}. Thus $a_1=\mu_1+2, a_2=\mu_2+2, b_1=\nu_1+1,b_2=\nu_2+1,b_3=\nu_3+0$ and $m(G)=4$. If $m=2$, then we have
$h(G)=(-q)^4\cdot(1,3,3,3,4).$

\bibliographystyle{plain}

\end{document}